\documentclass[reqno]{amsart}

\usepackage{tikz-cd}

\usepackage{fullpage,dsfont}

\usepackage{hyperref} 

\usepackage{parskip}
\makeatletter 
\def\thm@space@setup{%
 \thm@preskip=\parskip \thm@postskip=0pt
}
\def\th@remark{%
  \thm@headfont{\itshape}%
  \normalfont 
  \thm@preskip\parskip \thm@postskip=0pt
}
\makeatother

\usepackage[nobysame,alphabetic,initials,msc-links]{amsrefs}

\usepackage[final]{pdfpages}

\usepackage{multirow}

\usepackage{array}
\usepackage{kbordermatrix}

\DefineSimpleKey{bib}{how}
\DefineSimpleKey{bib}{mrclass}
\DefineSimpleKey{bib}{mrnumber}
\DefineSimpleKey{bib}{fjournal}
\DefineSimpleKey{bib}{mrreviewer}

\renewcommand{\PrintDOI}[1]{%
  \href{http://dx.doi.org/#1}{{\tt DOI:#1}}%
}
\renewcommand{\eprint}[1]{#1}
\BibSpec{book}{%
    +{}  {\PrintPrimary}                {transition}
    +{.} { \PrintDate}                  {date}
    +{.} { \textit}                     {title}
    +{.} { }                            {part}
    +{:} { \textit}                     {subtitle}
    +{,} { \PrintEdition}               {edition}
    +{}  { \PrintEditorsB}              {editor}
    +{,} { \PrintTranslatorsC}          {translator}
    +{,} { \PrintContributions}         {contribution}
    +{,} { }                            {series}
    +{,} { \voltext}                    {volume}
    +{,} { }                            {publisher}
    +{,} { }                            {organization}
    +{,} { }                            {address}
    +{,} { }                            {status}
    +{,} { \PrintDOI}                   {doi}
    +{,} { \PrintISBNs}                 {isbn}
    +{}  { \parenthesize}               {language}
    +{}  { \PrintTranslation}           {translation}
    +{;} { \PrintReprint}               {reprint}
    +{.} { }                            {note}
    +{.} {}                             {transition}
    +{}  {\SentenceSpace \PrintReviews} {review}
}
\BibSpec{article}{%
    +{}  {\PrintAuthors}                {author}
    +{,} { \textit}                     {title}
    +{.} { }                            {part}
    +{:} { \textit}                     {subtitle}
    +{,} { \PrintContributions}         {contribution}
    +{.} { \PrintPartials}              {partial}
    +{,} { }                            {journal}
    +{}  { \textbf}                     {volume}
    +{}  { \PrintDatePV}                {date}
    +{,} { \issuetext}                  {number}
    +{,} { \eprintpages}                {pages}
    +{,} { }                            {status}
    +{,} { \PrintDOI}                   {doi}
    +{,} { \eprint}        {eprint}
    +{}  { \parenthesize}               {language}
    +{}  { \PrintTranslation}           {translation}
    +{;} { \PrintReprint}               {reprint}
    +{.} { }                            {note}
    +{.} {}                             {transition}
    +{}  {\SentenceSpace \PrintReviews} {review}
}
\BibSpec{collection.article}{%
    +{}  {\PrintAuthors}                {author}
    +{,} { \textit}                     {title}
    +{.} { }                            {part}
    +{:} { \textit}                     {subtitle}
    +{,} { \PrintContributions}         {contribution}
    +{,} { \PrintConference}            {conference}
    +{}  {\PrintBook}                   {book}
    +{,} { }                            {booktitle}
    +{,} { \PrintDateB}                 {date}
    +{,} { pp.~}                        {pages}
    +{,} { }                            {publisher}
    +{,} { }                            {organization}
    +{,} { }                            {address}
    +{,} { }                            {status}
    +{,} { \PrintDOI}                   {doi}
    +{,} { \eprint}        {eprint}
    +{}  { \parenthesize}               {language}
    +{}  { \PrintTranslation}           {translation}
    +{;} { \PrintReprint}               {reprint}
    +{.} { }                            {note}
    +{.} {}                             {transition}
    +{}  {\SentenceSpace \PrintReviews} {review}
}
\BibSpec{misc}{%
  +{}{\PrintAuthors}  {author}
  +{,}{ \textit}      {title}
  +{.}{ }             {how}
  +{}{ \parenthesize} {date}
  +{,} { available at \eprint}        {eprint}
  +{,}{ available at \url}{url}
  +{,}{ }             {note}
  +{.}{}              {transition}
}
\usepackage{amssymb, amsfonts, amsxtra, amsmath}
\usepackage{mathrsfs}
\usepackage{mathdots}
\usepackage{wasysym}
\usepackage[all]{xy}
\usepackage{bbm}
\usepackage{calc}
\usepackage{accents}

\usepackage{stackengine}

\numberwithin{equation}{section}

\newtheorem{TheoIntro}{Theorem}

\newtheorem{Theorem}{Theorem}[section]
\newtheorem*{Theorem*}{Theorem}
\newtheorem{Def}[Theorem]{Definition}
\newtheorem*{Def*}{Definition}
\newtheorem{Lem}[Theorem]{Lemma}
\newtheorem{Prop}[Theorem]{Proposition}
\newtheorem{Cor}[Theorem]{Corollary}

\newtheorem{Rem}[Theorem]{Remark}

\newtheorem{Exa}[Theorem]{Example}

\newcommand\bp{\begin{proof}}
\newcommand\ep{\end{proof}}

\mathchardef\mhyph="2D

\DeclareMathOperator{\diag}{\mathrm{diag}}
\DeclareMathOperator{\ad}{\mathrm{ad}}

\DeclareMathOperator{\Ad}{\mathrm{Ad}}

\DeclareMathOperator{\id}{\mathrm{id}}

\DeclareMathOperator{\Tr}{\mathrm{Tr}}

\DeclareMathOperator{\braid}{\mathrm{br}}

\DeclareMathOperator{\Ker}{\mathrm{Ker}}

\DeclareMathOperator{\Det}{\mathrm{Det}}
\DeclareMathOperator{\Sym}{\mathrm{Sym}}

\DeclareMathOperator{\Sign}{\mathrm{Sign}}

\newcommand{\op}{\mathrm{op}}
\newcommand{\wt}{\mathrm{wt}}

\newcommand{\msR}{\mathscr{R}}

\newcommand{\mfl}{\mathfrak{l}}

\newcommand{\mft}{\mathfrak{t}}
\newcommand{\mfu}{\mathfrak{u}}

\newcommand{\mfgl}{\mathfrak{gl}}
\newcommand{\mfX}{\mathfrak{X}}

\newcommand{\mcC}{\mathcal{C}}

\newcommand{\Hsp}{\mathcal{H}}

\newcommand{\mcO}{\mathcal{O}}

\newcommand{\mcS}{\mathcal{S}}

\newcommand{\mbr}{\mathbf{r}}

\newcommand{\brr}{[r]}
\newcommand{\brk}{[k]}

\newcommand{\brt}{[t]}

\newcommand{\brN}{[N]}

\newcommand{\C}{\mathbb{C}}

\newcommand{\R}{\mathbb{R}}
\newcommand{\T}{\mathbb{T}}

\newcommand{\Z}{\mathbb{Z}}

\newcommand{\sa}{\mathrm{sa}}

\newcommand{\br}{\mathrm{br}}

\newcommand{\RE}{\mathrm{RE}}

\title{Representation theory of the reflection equation algebra II: Theory of Shapes}

\author{Kenny De Commer}
\address{Vrije Universiteit Brussel}
\email{kenny.de.commer@vub.be}

\author{Stephen T. Moore}
\address{Institute of Mathematics, Polish Academy of Sciences}
\email{stm862@gmail.com}

\thanks{The work of K.DC.~ and S.T.M.~was supported by the FWO grants G025115N and G032919N. S.T.M. was additionally supported by Narodowe Centrum Nauki, grant number 2017/26/A/ST1/00189.}

\begin{document}
\maketitle

\begin{abstract}
We continue our study of the representations of the Reflection Equation Algebra (=REA) on Hilbert spaces, focusing again on the REA constructed from the $R$-matrix associated to the standard $q$-deformation of $GL(N,\C)$ for $0<q<1$. We consider the Poisson structure appearing as the classical limit of the $R$-matrix, and parametrize the symplectic leaves explicitly in terms of a type of matrix we call a \emph{shape matrix}. We then introduce a quantized version of the shape matrix for the REA, and show that each irreducible representation of the REA has a unique shape. 
\end{abstract}

\section*{Introduction}

Let  $0<q<1$. We let $\hat{R} \in M_N(\C) \otimes M_N(\C)$ be the \emph{braid operator}
\begin{equation}\label{EqBraidOp}
\hat{R} =  \sum_{ij} q^{-\delta_{ij}}e_{ji}\otimes e_{ij} + (q^{-1} -q)\sum_{i<j} e_{jj} \otimes e_{ii}. 
\end{equation}
We recall that it is self-adjoint, $\hat{R}^* = \hat{R}$, and satisfies the \emph{braid relation}
\[
\hat{R}_{12}\hat{R}_{23}\hat{R}_{12} = \hat{R}_{23} \hat{R}_{12}\hat{R}_{23}.
\] 

The \emph{Reflection Equation $*$-Algebra} ($*$-\emph{REA}) is then the universal unital complex $*$-algebra $\mcO_{q}(H_N(\C))$ generated by the matrix entries of a selfadjoint matrix of generators
\[
Z = \sum_{ij} e_{ij} \otimes Z_{ij} \in M_N(\C)\otimes \mcO_q^{\br}(H(N)) = M_N(\mcO_q^{\br}(H(N))),\qquad Z^* = Z,
\] 
with universal relations given by the \emph{reflection equation}
\begin{equation}\label{EqUnivRelZ}
\hat{R}_{12}Z_{23} \hat{R}_{12}Z_{23} = Z_{23} \hat{R}_{12}Z_{23}\hat{R}_{12}. 
\end{equation}

The $*$-algebra $\mcO_q(H(N))$ can be viewed as a $q$-\emph{deformation} of the $*$-algebra of polynomial functions on $H(N)$, the space of $N\times N$ self-adjoint matrices.

The current paper is a continuation of \cite{DCMo24a}, which began a comprehensive study of the representation theory of $\mcO_q(H(N))$. By \emph{representation}, we mean a unital $*$-homomorphism into the $*$-algebra of bounded operators on a Hilbert space $\Hsp$,
\[
\pi:\mcO_q(H(N))\rightarrow B(\Hsp).\]
A key result of \cite{DCMo24a} was to show that $\mcO_q(H(N))$ is \emph{type $I$}, so its representations are in principle classifiable. The main goal of the current paper is to introduce the notion of \emph{shape} for the reflection equation algebra. This is a key property of $\mcO_q(H(N))$-representations, and will be used in \cite{Mo24} to obtain a classification of the irreducible representations of $\mcO_q(H(N))$.

To motivate our definition of shapes for the reflection equation algebra, we first consider the \emph{quasi-classical limit}. Under the limit $q\rightarrow 1$, $\mcO_q(H(N))$ reduces to the $*$-algebra of complex valued polynomial functions on $H(N)$. From the $R$-matrix, we obtain a \emph{Poisson structure}, known as the \emph{Semenov-Tian-Shansky (STS) Poisson bracket} \cite{Mud06}. With this bracket, $H(N)$ can be viewed as a Poisson manifold, and will have a foliation of symplectic leaves. From the philosophy of Kirillov's orbit method, there should be a close correspondence between the symplectic leaves of $H(N)$ and the irreducible representations of $\mcO_q(H(N))$. By way of motivation, we hence formulate the following result, that is based on well-known ideas. 

\begin{TheoIntro}\label{TheoMo}
Denote by $U(N)$ the group of $N\times N$ unitary matrices, and by $T(N)$ the space of $N\times N$ upper triangular matrices with strictly positive diagonal entries. Then the symplectic leaves of $H(N)$ with respect to the STS Poisson bracket are labelled by the intersections of $U(N)$-orbits and $T(N)$-orbits in $H(N)$ under the (right) action\footnote{Note that this can be viewed a as a twisted right conjugation $\Ad_x(z) = \theta(x)^{-1}zx$, where $\theta(x) = (x^*)^{-1}$.}
\[
\Ad_x: H(N) \rightarrow H(N), \quad z \mapsto xzx^*,\qquad x\in GL(N,\C). 
\] 
\end{TheoIntro} 

We provide a more concrete parametrisation of the above symplectic leaves using \emph{spectral weights} and \emph{shapes}. 

First, observe that an $\Ad_{U(N)}$-orbit is determined by the list of eigenvalues of any matrix in this orbit, i.e.\
\[
H(N)/U(N) \cong \Lambda := \R^N/\Sym(N).
\]
We may view this list of eigenvalues as the \emph{weight} of a self-adjoint matrix. 

The \emph{shape} of a selfadjoint matrix is then obtained through an explicit parametrization of the orbits of $\Ad_{T(N)}$ by means of \emph{shape matrices}.

\begin{Def*}
A \emph{(self-adjoint) shape matrix} is a self-adjoint $N\times N$ matrix, with all entries either zero or unimodular, and with all eigenvalues in $\{-1,0,1\}$.
\end{Def*}
We will show then that an intersection of a $U(N)$-orbit and a $T(N)$-orbit is uniquely determined by a compatible pair of a weight and a shape matrix, the compatibility being that they have the same \emph{signature}. 

Considering the connection between symplectic leaves of $H(N)$ and irreducible representations of $\mcO_q(H(N))$, it is natural to expect a \emph{quantization} of $\Lambda$ to correspond to some form of central or highest weights. In fact, we have already seen this appear in terms of \emph{central $*$-characters} in \cite{DCMo24a}. Our focus then should be on obtaining an appropriate quantization of shape matrices. To achieve this, we consider certain elements of $\mcO_q(H(N))$ obtained as $q$-deformations of sub-determinants labelled by a shape matrix. Under this correspondence, the quantized shape matrix associated to a representation $\pi$ simultaneously describes a number of elements in $\Ker(\pi)$, and a family of operators that $q$-commute with all other elements of $\pi\left(\mcO_q(H(N))\right)$. Our second main result can then be stated as follows:
\begin{TheoIntro}
Every irreducible representation of $\mcO_q(H(N))$ has a unique shape.
\end{TheoIntro}
The analogue of Theorem \ref{TheoMo} in the quantum setting will be considered then in \cite{Mo24}.

The precise contents of this paper are as follows. In Section \ref{section quasi-classical limit}, we consider the quasi-classical limit of $\mcO_q(H(N))$ as a Poisson manifold, and classify its symplectic leaves. In Section \ref{section: Quantum Minors and Laplace Expansions}, we consider a family of elements of $\mcO_q(H(N))$ known as \emph{quantum minors}, and prove some of the properties concerning quantum minors that we will need for our results. In Section \ref{section: shape}, we define our quantized notion of shape, and prove the uniqueness of shapes for $\mcO_q(H(N))$-representations.

We refer to \cite{DCMo24a} for detailed background and conventions. We end this introduction with some common notations that we will make use of.

For a given set $T$, we denote by $\mathscr{P}(T)$ the powerset of $T$. When $k \in \Z_{\geq 0}$, we denote $\brk$ for the set $\{1,2,\ldots k\}$. We denote $\binom{\brN}{k} \subseteq \mathscr{P}(\brN)$ for the set of subsets $I \subseteq \brN$ with $|I| = k$. We view $I \in \binom{\brN}{k}$ as a totally ordered set, and write
\[
I = \{i_1<\ldots < i_k\}
\]
to indicate this. We view $\binom{\brN}{k}$ itself as a partially ordered set in two ways: 
\begin{itemize}
\item as a totally ordered set by the  lexicographic order, denoted $\leq$, so 
\begin{equation}\label{EqOrdTot}
I < J \qquad \textrm{if and only if}\qquad i_1\leq j_1,\ldots, i_{p-1}\leq j_{p-1},i_p < j_p\textrm{ for some }1\leq p \leq k;
\end{equation}
\item as a partially ordered set with order $\preceq$ given by  
\begin{equation}\label{EqOrdPart}
I \preceq J \qquad \textrm{if and only if}\qquad i_p \leq j_p\textrm{ for all }1\leq p \leq N.
\end{equation}
\end{itemize}
We extend these orders lexicographically to orders on $\binom{\brN}{k}\times \binom{\brN}{k}$, so for example
\[
(J,I) \preceq (J',I') \qquad \iff \qquad J \prec J' \textrm{ or }(J = J'\textrm{ and }I\preceq I'). 
\]

We further use the following notation:  For $I\in \binom{\brN}{k}$ and $K\in \binom{\brk}{l}$ we write
\begin{equation}\label{EqSetUnderUpper}
I_K = \{i_{p}\mid p \in K\},\qquad I^K = I \setminus I_K.
\end{equation}
We denote by $I\Delta J$ the \emph{symmetric difference}, i.e.
\[
I\Delta J = (I\cup J)\setminus(I\cap J).\]
For $I \in \binom{\brN}{k}$, we define the \emph{weight} of $I$ to be the sum of its elements, 
\[
\wt(I) = \sum_{r=1}^k i_r.
\] 
For $X \in M_N(\C)$ and $I,J \in \binom{[N]}{k}$, we denote 
\[
X[I,J] \in M_k(\C)
\]
for the submatrix with $X[I,J]_{r,s} = X_{i_r,j_s}$. The associated minor is written as 
\[
X_{I,J} = \Det(X[I,J]). 
\]
Finally, for a bijection $\sigma: A \rightarrow B$ between finite totally ordered sets, we denote
\begin{equation}\label{EqLengthSym}
	l(\sigma) = |\{(i,j) \in A\times A \mid i<j,\sigma(i)>\sigma(j)\}|
\end{equation}
to be the number of inversions in $\sigma$.

\section{The quasi-classical limit}\label{section quasi-classical limit}

We can interpret \eqref{EqBraidOp} with $\hat{R} \in M_N(\C[[\hbar]]) \otimes_{\C[[\hbar]]} M_N(\C[[\hbar]])$ by putting $q = e^{\hbar}$. Here we view $\hbar$ as a real formal variable, i.e.\ we have the $*$-structure $\hbar^* = \hbar$. The defining relation \eqref{EqUnivRelZ}  gives us the topological $\C[[\hbar]]$-$*$-algebra $\mcO_{\hbar}^{\braid}(H(N))$, whose generators we denote by $Z_{\hbar,ij}$. Denoting by $Z_{ij} \in \mcO(H(N))$ the coordinate functions on $H(N)$, and noting that we can identify $Z_{ij} = Z_{\hbar,ij} \!\!\mod \hbar$, we have that 
\begin{equation}
[Z_{\hbar,ij},Z_{\hbar,kl}] \mod \hbar^2 = i\hbar \{Z_{ij},Z_{kl}\}_{\RE},
\end{equation}
where $\{-,-\}_{\RE}$ is the unique Poisson bracket on the polynomial algebra $\mcO(H(N))$ determined by 
\begin{equation}
\sum_{ijkl} e_{ij}\otimes e_{kl} \otimes \{Z_{ij},Z_{kl}\}_{\RE} = -i \lbrack (Z_{13}r_{12})^* + Z_{13}r_{12},Z_{23}],
\end{equation}
with $r$ the \emph{classical $r$-matrix}
\begin{equation}
r = \sum_{i}e_{ii}\otimes e_{ii} + 2 \sum_{i<j}e_{ij}\otimes e_{ji}.
\end{equation}
Identifying $H(N)$ with its tangent space at any point, the Poisson bivector $\mfX^{\RE}$ associated to $\{-,-\}_{\RE}$ is given
at $z \in H(N)$ by 
\[
i\mfX_z^{\mathrm{RE}} = r_{21}(z\otimes z) - (z\otimes z)r+ (z\otimes 1)r(1\otimes z) -(1\otimes z)r_{21}(z\otimes 1) .
\]

We view $(H(N),\{-,-\}_{\RE})$ as the quasi-classical limit of $\mcO_q(H(N))$. Based on the philosophy of Kirillov's orbit method, one expects to see a close connection between the representation theory of $O_q(H(N))$ on the one hand, and the structure of the symplectic leaves in $H(N)$ for $\{-,-\}_{\mathrm{RE}}$ on the other. In this section, we present some details on the latter. 

Consider $H(N)$ with the adjoint action by $GL(N,\C)$ given by 
\begin{equation}\label{EqAdjAct}
\Ad_x(z) = x^*zx.
\end{equation}
Let $U(N) \subseteq GL(N,\C)$ be the group of unitary matrices and $T(N)\subseteq GL(N,C)$ the group of upper triangular matrices with strictly positive diagonal entries.

\begin{Def}
We define  $\widehat{H}(N) = \{C\} \subseteq \mathscr{P}(H(N))$ as the non-empty subsets arising as the intersection of an $\Ad_{U(N)}$-orbit and an $\Ad_{T(N)}$-orbit under the adjoint action \eqref{EqAdjAct}.
\end{Def}
We thus obtain a partition 
\[
H(N) \twoheadrightarrow \widehat{H}(N),\qquad z \mapsto [z] := \Ad_{T(N)}(z) \cap \Ad_{U(N)}(z) \subseteq H(N).
\]

\begin{Prop}\label{PropSympLeaves}
The $C \in \widehat{H}(N)$ are connected, and form the symplectic leaves for $\{-,-\}_{\mathrm{RE}}$. 
\end{Prop}
\begin{proof}
Let us first show that the $C \in \widehat{H}(N)$ are indeed connected. Take $z\in H(N)$. Then we have a continuous surjection 
\[ 
X_z = \{(t,u) \in 
T(N) \times U(N)\mid  t^*zt = uzu^*\} \twoheadrightarrow [z],\qquad (t,u) \mapsto uzu^*. 
\] 
But as the multiplication map gives a homeomorphism $T(N) \times U(N) \rightarrow GL(N,\C)$, we see that 
\[
X_z\cong \{x \in GL(N,\C) \mid x^*zx = z\}.
\] 
As the latter spaces are all homeomorphic for $z$'s of the same signature, we may as well assume that $z$ is of the form $(\underbrace{+,\ldots,+}_{N_+},\underbrace{-,\ldots,-}_{N_-},\underbrace{0,\ldots,0}_{N_0})$, in which case $X_z \cong \begin{pmatrix} U(N_+,N_-) & M_{N_++N_-,N_0}(\C)\\ 0& GL(N_0,\C)\end{pmatrix}$, a connected topological space. 

To see that the $C \in \widehat{H}(N)$ are symplectic submanifolds of $H(N)$, we note first that the $U(N)$-orbits of any $z\in H(N)$ are Poisson submanifolds of $H(N)$. Indeed,  writing $r = \sum_{ij} a_{ij} e_{ij}\otimes e_{ji}$, a small computation reveals that
\[
2\mfX^{\RE} = \sum_{ij} a_{ij} \left((ize_{ij} -ie_{ji}z)\otimes [e_{ij}-e_{ji},z] - (ze_{ij} +e_{ji}z)\otimes [ie_{ij}+ie_{ji},z]\right),
\]
where $[e_{ij}-e_{ji},z],[ie_{ij}+ie_{ji},z] \in \ad_{\mfu(N)}(z)$. 

Recall now that, with $\mft(N)$ the Lie algebra of $T(N)$, the triple $(\mfgl(N,\C),\mfu(N),\mft(N))$ forms a (real) Manin triple with respect to the inner product $\langle X,Y \rangle = \mathrm{Im}(\Tr(XY))$ on $\mfgl(N,\C)$. The resulting Lie bialgebra structure on $\mfgl(N,\C)$ extends to a Poisson-Lie group structure on $GL(N,\C)$, and it is standard to verify that the adjoint action on $H(N)$ is then a Poisson action. 

Let $O_z$ be the $U(N)$-orbit of $z \in H(N)$. Then we obtain by restriction of the adjoint action a Poisson action of $U(N)$ on $O_z$. We may here of course assume that $z = \diag(z_1,\ldots,z_N)$ with signature $\diag(\underbrace{1,\ldots,1}_{N_+},\underbrace{-1,\ldots,-1}_{N_-},\underbrace{0,\ldots,0}_{N_0})$. It is then easily computed, following the terminology of \cite{Dri93}, that the Lagrangian subalgebra $\mfl$ of $\mfgl(N,\C)$ corresponding to the Poisson $U(N)$-homogeneous space $O_z$ is given by the Lie algebra of $L = \begin{pmatrix} U(N_+,N_-) & M_{N_++N_-,N_0}(\C)\\ 0& U(N_0) \end{pmatrix}$. In particular, writing 
\[
K = L \cap U(N) = U(N_+) \times U(N_-) \times U(N_0),
\]
we obtain a Poisson embedding 
\[
\iota_z: O_z \cong K\backslash U(N) \rightarrow L \backslash GL(N,\C). 
\] 
By \cite{Kar95} (see also \cite[Remark 2.10]{LY08}), the symplectic leaves of $O_z$ are given by connected components of the inverse image under $\iota_z$ of the $T(N)$-orbits in $L \backslash GL(N,\C)$. Writing $L_0 = \begin{pmatrix} U(N_+,N_-) & M_{N_++N_-,N_0}(\C)\\ 0& GL(N_0,\C)\end{pmatrix}$, we get a factorisation of the inclusion $O_z \subseteq H(N)$ through the maps
\[
O_z \hookrightarrow L \backslash GL(N,\C) \twoheadrightarrow L_0 \backslash GL(N,\C) \hookrightarrow H(N)
\]
where the middle map is $GL(N,\C)$-equivariant. This implies that the symplectic leaves of $O_z$ are the connected components of intersections with the $T(N)$-orbits in $H(N)$ under the adjoint action, finishing the proof of the proposition. 
\end{proof}

We want to obtain a concrete parametrisation of the elements of $\widehat{H}(N)$. Given how $\widehat{H}(N)$ is defined, this requires the following: 
\begin{enumerate}
\item Parametrise the $U(N)$-orbits.
\item Parametrise the $T(N)$-orbits. 
\item Determine when a $U(N)$-orbit and $T(N)$-orbit have non-empty intersection. 
\end{enumerate}

The first item is easily treated. 
\begin{Def} We define
\[
\Lambda = \R^N/\Sym(N),
\]
the multisets of $N$ real numbers. We call $\Lambda$ the space of \emph{spectral weights}. 
\end{Def}
By the spectral theorem, we have an isomorphism 
\[
H(N)/U(N)\cong \Lambda,\qquad \Ad_{U(N)}(z) \mapsto \lambda(z),
\]
with $\lambda(z)$ the multiset of eigenvalues of $z$. In the following, we will identify $\Lambda$ with the set of diagonal matrices in $H(N)$ whose values along the diagonal are non-decreasing. 

To parametrise the $T(N)$-orbits, we need some more preparation. Define 
\[
\T = \{u \in \C\mid |u|= 1\}.
\]

\begin{Def}\label{DefShape}
A \emph{shape} is a couple $S=(\tau,u)$ with 
\begin{itemize}
\item $\tau$ an involution of $\{1,2,\ldots,N\}$ and 
\item $u: \{1,2,\ldots,N\} \rightarrow \{0\}\cup \T$ such that $u_i \neq 0$ if $i\neq \tau(i)$.
\end{itemize}
We write $P = \{p \in \brN \mid u_p \neq 0\}$. We call $P$ the \emph{support set} and $M=|P|$ the \emph{rank} of $S$. 

We call a shape $S$ \emph{self-adjoint} if $u_{\tau(i)} = \overline{u_i}$. 
\end{Def}
We denote 
\begin{equation}
\mcS = \{\textrm{ shapes }\},\qquad \mcS_{\sa}= \{\textrm{ self-adjoint shapes }\}.
\end{equation}

We can encode a shape $S$ by its matrix $S \in M_N(\{0\}\cup \T)$ with $Se_i =u_i e_{\tau(i)}$. Then $S$ is self-adjoint if and only if $S$ is self-adjoint as a matrix, and the rank of $S$ is its rank as a matrix. In the following, we will always use this identification of shapes with matrices implicitly. It is elementary to verify that the resulting class of self-adjoint shape matrices indeed coincides with the class described in the introduction.

Let $S = (\tau,u)$ be a shape with support set $P$ as above. If $P = \{p_1<p_2<\ldots<p_M\}$ we further write 
\[
P_{\brk} = \{p_1,p_2,\ldots,p_k\},\qquad \tau(P_{\brk}) = \{\tau(p_1),\tau(p_2),\ldots,\tau(p_k)\}.
\]
For $Q \subseteq \{1,2,\ldots,N\}$, we write 
\[
u_Q := \prod_{q\in Q} u_q.
\]
Recall from \eqref{EqOrdTot} and \eqref{EqOrdPart} the orders $\preceq$ and $\leq$ on $\binom{\brN}{k}$, and the function $l$ defined by \eqref{EqLengthSym}. 

\begin{Def}
A matrix $z \in H(N)$ is said to be of shape $S = (\tau,u)$ if the rank of $S$ equals the rank of $z$ and, for all $k \leq M = |S|$,
\begin{itemize}
\item for $I,J \in \binom{[M]}{k}$ the minors $z_{I,J}$ vanish for $(J,I) < (P_{[k]},\tau(P_{[k]}))$, and
\item $(-1)^{l(\tau_{|P_{[k]}})}\overline{u_{P_{[k]}}} z_{\tau(P_{[k]}),P_{[k]}}>0$. 
\end{itemize}
\end{Def}

\begin{Theorem}
Each $z \in H(N)$ has a unique shape $S = S(z)$, and $S(z)$ is self-adjoint. Moreover, the map 
\[
\mcS_{\sa} \rightarrow H(N)/T(N),\qquad S \mapsto \Ad_{T(N)}(S)
\]
is bijective, with inverse map $\Ad_{T(N)}(z) \mapsto S(z)$. 
\end{Theorem}

\begin{proof}
We first show that all elements in the same $T(N)$-orbit share the same shape. For this, it is sufficient to see that associated shapes are invariant under conjugation by diagonal matrices and matrices $I + \lambda e_{r,r+1}$ for $\lambda \in \C$, as these generate $T(N)$. Invariance under diagonal matrices is clear, so fix $t= I + \lambda e_{r,r+1}$. It is sufficient to prove the following claim: For $x \in M_N(\C)$ and $(J,I) \in \binom{\brN}{k} \times \binom{\brN}{k}$ the $\leq$-first element with $0\neq x_{I,J}$, we have that $(J,I)$ is also the  $\leq$-first element with $0\neq x_{I,J}'$ for $x' = xt$ or $x' = t^*x$, and moreover $x_{I,J}' = x_{I,J}$. But, for general $K,L$, we easily see that 
\[
(xt)_{K,L} = 
\left\{\begin{array}{ll} 
x_{K,L} &\textrm{if }r+1\notin L \textrm{ or }\{r,r+1\}\subseteq L, \\
x_{K,L} + \lambda x_{K,(L\setminus \{r+1\})\cup \{r\}} &\textrm{if }r\notin L,r+1\in L,
\end{array}\right.
\]
where we note that in the second case $(L\setminus \{r+1\})\cup \{r\} <L$. Similarly, 
\[
(t^*x)_{K,L} = 
\left\{\begin{array}{ll} 
x_{K,L} &\textrm{if }r+1\notin K\textrm{ or }\{r,r+1\}\subseteq K, \\
x_{K,L} + \lambda x_{(K\setminus \{r+1\})\cup \{r\},L} &\textrm{if }r\notin K,r+1\in K.
\end{array}\right.
\]
This proves the claim. 

Note now that by \cite{Els79}, or more conveniently \cite[Theorem 7]{GG82} for the specific case of self-adjoint matrices, there exists for any $z \in H(N)$ a unique\footnote{In the above references the gauge factor $u$ can be made real as $t$ is allowed to have complex values on the diagonal.} self-adjoint shape $S\in \mcS_{\sa}$ such that
\begin{equation}\label{EqEls}
z = t^*St,\qquad t \in T(N).
\end{equation} 

To finish the proof, it is sufficient to show that any $S\in \mcS_{\sa}$, interpreted as an element of $H(N)$, has $S$ as its unique associated shape. This is however clear by direct inspection.
\end{proof}

Finally, let us address when a $T(N)$-orbit and a $U(N)$-orbit in $H(N)$ have non-empty intersection. Let $\Sign$ be the multiset of $N$ elements chosen from $\{+,-,0\}$. For $X$ a selfadjoint matrix, let $\Sign(X)$ be the multiset of signs of its eigenvalues. Define 
\[
\mcC = \mcS_{\sa}\underset{\mathrm{Sign}}{\times} \Lambda = \{(S,\lambda) \mid \mathrm{Sign}(S) = \mathrm{Sign}(\lambda)\}. 
\]

\begin{Theorem}\label{TheoRhoBij}
The map
\[
\rho: \widehat{H}(N) \rightarrow \mcC,\quad [z]\mapsto (S(z),\lambda(z)).
\]
is well-defined and bijective.
\end{Theorem} 
\begin{proof}
As $S(z)$ lies in the $T(N)$-orbit of $z\in H(N)$ by the (proof of the) previous theorem, we have that $\mathrm{Sign}(z) = \mathrm{Sign}(S(z)) = \mathrm{Sign}(\lambda(z))$, so the map $\rho$ is well-defined. 

 To see that $\rho$ is surjective, note first that for $\lambda_1,\lambda_2>0$ and $u\in \T$, we have that
\begin{equation}\label{EqEigGen}
\begin{pmatrix} \sqrt{\lambda_1\lambda_2} & \frac{u}{2}(\lambda_1-\lambda_2)  \\ 0& 1\end{pmatrix}^* \begin{pmatrix} 0 & u \\ \overline{u} & 0 \end{pmatrix} \begin{pmatrix} \sqrt{\lambda_1\lambda_2} & \frac{u}{2}(\lambda_1-\lambda_2)  \\ 0& 1\end{pmatrix}
\end{equation}
is a matrix with eigenvalues $\lambda_1,-\lambda_2$. Hence, by looking at the $1$-by-$1$ or $2$-by-$2$-blocks appearing in a shape $S$, we can, by conjugation by an upper triangular matrix $t \in T(N)$ with the same block arrangement, arrange for $t^*St$ to have any list of eigenvalues compatible with the signature. More precisely, $t^*st$ will then consist of diagonal elements $(\alpha_i)_{ii}$ for $\alpha_i \in \R$ when $i = \tau(i)$, and $2$-by-$2$-blocks
\renewcommand{\kbldelim}{(}
\renewcommand{\kbrdelim}{)}
\[
\kbordermatrix{
& i & \tau(i) \\
i & 0 & \beta_i \\ 
\tau(i) & \overline{\beta_i} & \delta_i
}
\]
for $\delta_i \in \R$ and $\beta_i \in \C^{\times}$ when $i < \tau(i)$. Call a self-adjoint matrix of this form an \emph{enhanced self-adjoint shape}. 

It then follows from the above that $\rho$ is indeed surjective, and that in fact any element of $\widehat{H}(N)$ contains an enhanced self-adjoint shape $E$. If now $E,E'$ are enhanced self-adjoint shapes with the same shape $S_E = S_{E'}$ and the same list of eigenvalues $\lambda(E) = \lambda(E')$, it is then again immediate, by restricting to $1$-by-$1$- and $2$-by-$2$-blocks and using \eqref{EqEigGen} for the latter, that $[E] = [E']$, proving injectivity of $\rho$.
\end{proof}

We can summarize the above by saying that there is a bijection of pullback diagrams
\begin{equation}\label{EqFundSquare}
\xymatrix{&\widehat{H}(N)  \ar[rd]\ar[ld]&     &   &  &\mcC \ar[rd]\ar[ld]& \\ 
H(N)/T(N) \ar[rd]& & H(N)/U(N)\ar[ld]  & \cong & \mcS \ar[rd]&& \Lambda \ar[ld]\\ 
& H(N)/GL(N,\C) & && &\mathrm{Sign}&}
\end{equation}

\section{Quantum Minors and Laplace Expansions}\label{section: Quantum Minors and Laplace Expansions}

We recall some terminology and results from \cite{DCMo24a}. 

\subsection{Quantum minors in quantum matrices}

The \emph{quantum Grassmann algebra} $\wedge_q(\C^N)$ is generated by elements $e_i$ for $1\leq i\leq N$, with product $\wedge_q$ satisfying the universal relations 
\[
e_i \wedge_q e_j = -qe_i \wedge_q e_j,\qquad 1\leq i\leq j \leq N. 
\]
It has the basis
\begin{equation}\label{EqGrassBasis}
e_I = e_{i_1}\wedge_q \ldots \wedge_q e_{i_k},\qquad I \in \binom{[N]}{k}, \qquad I = \{i_1<\ldots<i_k\}.
\end{equation}
With $\mcO_q(M_N(\C))$ the unique bialgebra generated by the universal matrix $X = (X_{ij})_{1\leq i,j\leq N}$ with 
\[
\hat{R}_{12}X_{13}X_{23} =X_{13}X_{23}\hat{R}_{12},\qquad (\id\otimes \Delta)X = X_{12}X_{13},
\]
we have that $\wedge_q(\C^N)$ carries a unique right $\mcO_q(M_N(\C))$-comodule structure $\wedge \delta$ such that
\begin{equation}\label{EqComodAlgMapqExt}
	\wedge \delta: \wedge_q(\C^N) \rightarrow \wedge_q(\C^N)\otimes \mcO_q(M_N(\C)),\qquad (\wedge \delta)e_i = \sum_{j} e_j\otimes X_{ji}. 
\end{equation}
We then write the  \emph{quantum minors} $X_{IJ}\in\mcO_q(U(N))$ as the matrix coefficients with respect to the basis \eqref{EqGrassBasis},
\[
(\wedge \delta)(e_I) = \sum_{J\in \binom{[N]}{k}} e_J \otimes X_{JI},\qquad I \in \binom{[N]}{k}. 
\]

Inside $\mcO_q(M_N(\C))$, the quantum minors satisfy a braided commutativity relation, given in general by 
\begin{equation}\label{EqBraidComMInvAlt}
	ab = \mbr^{-1}(a_{(1)},b_{(1)})b_{(2)}a_{(2)} \mbr(a_{(3)},b_{(3)}) = \mbr(b_{(1)},a_{(1)})b_{(2)}a_{(2)} \mbr^{-1}(b_{(3)},a_{(3)}),\qquad a,b\in \mcO_q(M_N(\C)). 
\end{equation}
Here $\mbr$ is the skew bicharacter uniquely determined by 
\begin{multline*}
(\id\otimes \id\otimes \mbr)(X_{13}X_{24}) = \Sigma\circ \hat{R},\quad \mbr(ab,c) = \mbr(a,c_{(1)})\mbr(b,c_{(2)}),\quad \mbr(a,bc) = \mbr(a_{(2)},b)\mbr(a_{(1)},c),\\ a,b,c\in \mcO_q(M_N(\C)), 
\end{multline*}
where $\Sigma$ is the flip map. For $I,J\in \binom{[N]}{k}$ and $I',J'\in \binom{[N]}{l}$, we then denote
\begin{equation}\label{EqCoeffsRBraid}
	\hat{R}^{IJ}_{I'J'} := \mbr(X_{JI},X_{I'J'}),\qquad (\hat{R}^{-1})^{IJ}_{I'J'} = \mbr^{-1}(X_{JI},X_{I',J'}).
\end{equation}

As in Lemma $2.7$ of \cite{DCMo24a}, the $\hat{R}^{IJ}_{I'J'}$ satisfy
\begin{equation}\label{EqCondR}
	\hat{R}^{IJ}_{I'J'} \neq 0 \quad \textrm{only if}\quad  J\preceq I,J'\preceq I'\quad \textrm{and}\quad J \setminus I = J'\setminus I',\;I\setminus J = I'\setminus J',
\end{equation}
and
\[
\hat{R}^{II}_{I'I'} = q^{- |I\cap I'|}. 
\] 
The values $(\hat{R}^{-1})_{I'J'}^{IJ}$ also satisfy the conditions of \eqref{EqCondR}, however in this case 
\[
(\hat{R}^{-1})^{II}_{I'I'} = q^{|I\cap I'|}.
\]

From \eqref{EqBraidComMInvAlt}, we have the following concrete form of the braided commutativity.
\begin{Prop}\label{PropBraidComm}
	For all $I,J \in \binom{\brN}{k}$ and $I',J'\in \binom{\brN}{l}$, one has
	\begin{eqnarray}
		X_{I',J'} X_{I,J} &=& \sum_{A,B,C,D} \hat{R}^{A,I}_{I',B} (\hat{R}^{-1})^{J,C}_{D,J'} X_{A,C} X_{B,D} \label{EqBraidComm1} \\
		& = &\sum_{A,B,C,D} (\hat{R}^{-1})^{B,I'}_{I,A} \hat{R}^{J',D}_{C,J} X_{A,C} X_{B,D}. \label{EqBraidComm2}
	\end{eqnarray}
\end{Prop} 

The following Laplace expansion is well-known \cite{PW91}: 
\begin{Prop}\label{PropLaplace}
	For $I,J\in \binom{\brN}{k}$ and $K,K'\in \binom{\brk}{l}$ the following identities hold in $\mcO_q(M_N(\C))$,
	\begin{eqnarray}
		\delta_{K,K'}X_{I,J} &=& \sum_{P \in \binom{\brk}{l}} (-q)^{\wt(P) -\wt(K)}X_{I_K,J_P}\; X_{I^{K'},J^P} \label{EqLaplaceRow}\\
		&=& \sum_{P \in \binom{\brk}{l}} (-q)^{\wt(P) -\wt(K)}X_{I_P,J_K}\; X_{I^P,J^{K'}}.  \label{EqLaplaceColumn}
	\end{eqnarray}
\end{Prop}

We will need the following generalisation, known as the `common-submatrix Laplace expansion'. It is attributed to Muir \cite{Mui60,Wil86,Wil15}. It follows immediately from the more general `quantum Muir law of extensible minors'  \cite[Theorem 3.8]{KL95}, applied to the Laplace expansion.

\begin{Prop}\label{PropMuir}
	Let $I,J \in \binom{\brN}{k}$ and $F,G\in \binom{[k]}{k-r}$. Then for any $K,K' \in \binom{\brr}{l}$ we have 
	\begin{eqnarray}\label{EqMuir}
		\delta_{K,K'}X_{I,J}X_{I_F,J_G} &=& \sum_{P \in \binom{\brr}{l}} (-q)^{\wt(P) -\wt(K)}X_{I_F \cup(I^F)_K,J_G\cup(J^G)_P}X_{I_F\cup (I^F)^{K'},J_G\cup (J^G)^P} \\
		&=& \sum_{P \in \binom{\brr}{l}} (-q)^{\wt(P) -\wt(K)}X_{I_F \cup(I^F)_P,J_G\cup(J^G)_K}X_{I_F\cup (I^F)^{P},J_G\cup (J^G)^{K'}} \label{EqMuir2}
	\end{eqnarray}
\end{Prop}

\subsection{Quantum minors in the reflection equation algebra}

There is also a convolution inverse $\mbr'$ of $\mbr$ when considered as a functional on $(\mcO_q(M_N(\C)),\Delta)\otimes(\mcO_q(M_N(\C)),\Delta^{\op})$. This allows us to define a new product $*$ on $\mcO_q(M_N(\C))$ via
\begin{equation}\label{EqBraidProd}
	f*g = \mbr(f_{(1)},g_{(2)})f_{(2)}g_{(3)} \mbr'(f_{(3)},g_{(1)}), 
\end{equation}
with conversely
\begin{equation}\label{EqRevBraid}
	fg = \mbr^{-1}(f_{(1)},g_{(1)})f_{(2)}*g_{(3)} \mbr(f_{(3)},g_{(2)}). 
\end{equation}
From \cite[Section 7.4]{Maj95}, \cite[Example 10.18]{KS97}, there exists a (unique) isomorphism of algebras
\begin{equation}\label{EqDefPhi}
	\Phi: \mcO_q(H(N)) \rightarrow (\mcO_q(M_N(\C)),*),\qquad Z \mapsto X. 
\end{equation}

In the case of quantum minors, \eqref{EqRevBraid} reads
\begin{equation}\label{EqRevBraidConc}
	X_{I,J}X_{I',J'} = \sum_{A,B,C,D} (\hat{R}^{-1})_{I',B}^{A,I} \hat{R}_{B,D}^{J,C} X_{A,C}*X_{D,J'}.
\end{equation}

\begin{Def}\label{DefQMinRefl}
The \emph{quantum minors} $Z_{IJ}\in\mcO_q(H(N))$ are defined by
\[
Z_{IJ}:=\Phi^{-1}(X_{IJ}).\]
\end{Def}

The analogue of the braided commutativity from \eqref{EqBraidComm2} then becomes 
\begin{equation}\label{EqGenCommRel}
	\sum_{K,L,L'} \left(\sum_{P'}\hat{R}_{JK}^{P'I'}\hat{R}_{P'L'}^{IL}\right)  Z_{K,L} Z_{L',J'} =   \sum_{K,L,L'}  \left(\sum_{P'}\hat{R}_{JK}^{P'L'}\hat{R}_{P'J'}^{IL}\right) Z_{I',L'} Z_{K,L}.
\end{equation}

We have an $\mcO_q(H(N))$-version of the Laplace expansion: 
\begin{Prop}
For $I,J \in \binom{\brN}{k}$ and $K\in \binom{\brk}{m}$ with $m\leq k \leq N$ we have
\begin{equation}\label{EqLaplExp1}
	Z_{I,J}= \sum_{P \in \binom{\brk}{m}}\sum_{S,T\in \binom{\brN}{m}}\sum_{S',T'\in \binom{\brN}{k-m}}(-q)^{\wt(P) -\wt(K)} (\hat{R}^{-1})^{S,I_K}_{I^K,T'}  \hat{R}^{J_P,T}_{T',S'} Z_{S,T} Z_{S',J^P}
\end{equation}
and
\begin{equation}\label{EqLaplExp2}
	Z_{I,J} =  \sum_{P \in \binom{\brk}{m}}\sum_{S,T\in \binom{\brN}{m}}\sum_{S',T'\in \binom{\brN}{k-m}}(-q)^{\wt(P) -\wt(K)} (\hat{R}^{-1})^{T,J_K}_{J^K,S'}  \hat{R}^{I_P,S}_{S',T'} Z_{I^P,T'} Z_{S,T}.
\end{equation}
\end{Prop}

We can now introduce the key identity that we will need for quantizing shapes. It is obtained from the  common-submatrix Laplace expansion of Proposition \ref{PropMuir} by applying equation \eqref{EqRevBraidConc}.

\begin{Prop}
	Let $I,J \in \binom{\brN}{k}$ and $F,G\in \binom{[k]}{k-r}$. Then for any $K,K' \in \binom{\brr}{l}$ we have 
	\begin{eqnarray} 
		& \hspace{-8cm}\delta_{K,K'}\sum_{S,T,H,L} (\widehat{R}^{-1})^{S,I}_{I_F,H} \widehat{R}^{J,T}_{H,L} Z_{S,T}Z_{L,J_G} \nonumber \\ 
		&= \sum_{A,B,C,D} \sum_{P \in \binom{\brr}{l}} (-q)^{\wt(P) -\wt(K)} (\widehat{R}^{-1})^{A,I_F\cup (I^F)_K}_{I_F\cup (I^F)^{K'},B} \widehat{R}^{J_G\cup (J^G)_P,C}_{B,D} Z_{A,C} Z_{D,J_G\cup (J^G)^P}  \label{EqMuirBr} \\
		&=  \sum_{A,B,C,D} \sum_{P \in \binom{\brr}{l}} (-q)^{\wt(P) -\wt(K)} (\widehat{R}^{-1})^{A,I_F\cup (I^F)_P}_{I_F\cup (I^F)^{P},B} \widehat{R}^{J_G\cup (J^G)_K,C}_{B,D} Z_{A,C} Z_{D,J_G\cup (J^G)^{K'}}. \label{EqMuirBr2} 
	\end{eqnarray}
\end{Prop}
We will refer to these expansions respectively as $L_l(I,J;F,G;K,K')$ and $L_r(I,J;F,G;K,K')$.

\section{Spectral shape}\label{section: shape}

We now introduce the quantized notion of \emph{shape} for the reflection equation algebra, and show that any irreducible representation of $\mcO_q(H(N))$ comes with an associated shape. 

\subsection{Rank and shape of representations}
We first recall from \cite{DCMo24a} the definition of rank for $\mcO_q(H(N))$ representations.
\begin{Def}
For $1\leq M\leq N$, the set of all $Z_{I,J}$ with $I,J \in \binom{\brN}{M}$ generates a 2-sided $*$-ideal of $\mcO_q(H(N))$ containing all $Z_{I'J'}\in\binom{\brN}{k}$ with $k\geq M$. We define a representation $\pi$ of $\mcO_q(H(N))$ to have rank $M$ if $I_{M+1}\subseteq\Ker(\pi)$ but $I_M\not\subseteq\Ker(\pi)$.
\end{Def}

\begin{Def}
Let $S$ be a shape of rank $M$. We define $I_S^{\preceq}$ to be the 2-sided $*$-ideal of $\mcO_q(H(N))$ generated by $I_{M+1}$ and, for each $k\leq M$, the $Z_{I,J}$ with $I,J \in \binom{\brN}{k}$ and $(J,I) \prec (P_{\brk},\tau(P_{\brk}))$. We denote
\[
\mcO_q^{S,\preceq}(H(N)) = \mcO_q(H(N))/I_S^{\preceq}.
\]   
We similarly define $I_S^{\leq}$ and $\mcO_q^{S,\leq}(H(N))$ by replacing $\preceq$ by $\leq$. 
\end{Def}
Clearly $I_S^{\preceq} \subseteq I_S^{\leq}$, hence we have a natural surjective $*$-homomorphism
\[
\mcO_q^{S,\preceq} \twoheadrightarrow \mcO_q^{S,\leq}.
\]

If $S$ is a shape of rank $M$, we will also use the notation 
\[
Z_{S,k} = Z_{\tau(P_{\brk}),P_{\brk}},\qquad k\leq M, 
\]
see Definition \ref{DefShape} and Definition \ref{DefQMinRefl}. 

\begin{Lem}\label{LemqCommpi}
In $\mcO_q^{S,\preceq}(H(N))$ we have the following commutation relations:
\begin{equation}\label{EqqCommQuot}
Z_{S,k} Z_{I,J} = q^{|I\cap P_{\brk}| + |I \cap \tau(P_{\brk})| - |J\cap P_{\brk}| - |J\cap \tau(P_{\brk})|} Z_{I,J}Z_{S,k}.
\end{equation}
\end{Lem}
\begin{proof}
This follows directly from \eqref{EqGenCommRel} upon taking $J = \tau(P_{\brk})$ and $I = P_{\brk}$.
\end{proof}

\begin{Cor}
For all $k,l\leq M$ we have $[Z_{S,k},Z_{S,l}] = [Z_{S,k},Z_{S,l}^*] =0$ in $\mcO_q^{S,\preceq}(H(N))$.
\end{Cor}

\begin{Def}
Let $S = (\tau,u)$ be a shape of rank $M$. We say that a representation $\pi$ of $\mcO_q(H(N))$ has \emph{shape $S$} if
\begin{itemize}
\item $I_S^{\leq} \subseteq \Ker(\pi)$,
\item $\Ker(\pi(Z_{S,k})) = 0$ for all $k\leq M$,
\item $(-1)^{l(\tau_{\mid P_{[k]}})}\overline{u_{P_{[k]}}}\pi(Z_{S,k})$ is a positive operator for all $k\leq M$.
\end{itemize}
\end{Def}

The last condition says that the unimodular part of $\pi(Z_{S,k})$ is a scalar equal to the determinant of the matrix $(S_{ij})_{i\in \tau(P_{[k]}),j\in P_{[k]}}$.

The following theorem establishes that an irreducible representation of $\mcO_q(H(N))$ has a unique shape.

\begin{Theorem}\label{TheoShapeIrr}
Each irreducible representation $\pi$ of $\mcO_q(H(N))$ has a unique shape $S$. Moreover, if $\pi$ has rank $M$, then $S$ has rank $M$. 
\end{Theorem}
\begin{proof}
It is immediate from the definitions that if $\pi$ is an irreducible representation with associated shape $S$, then the rank of $\pi$ is the rank of $S$. 

Fix now an irreducible representation $\pi$ of rank $M$. For each $k\leq M$ we define $J[k],I[k] \in \binom{\brN}{k}$ such that $(J[k],I[k])$ is the $\leq$-first couple $(J,I)$ for which $\pi(Z_{I,J})\neq 0$. Similar to Lemma \ref{LemqCommpi}, it is clear that $\pi(Z_{I[k],J[k]})$ are normal elements $q$-commuting with the $\pi(Z_{ij})$. Hence, by irreducibility, each $\pi(Z_{I[k],J[k]})$ will have zero kernel and purely discrete spectrum with constant unimodular part in its polar decomposition. The theorem will now follow if we can show the following three properties: 
\begin{enumerate}
\item\label{Prop1} $J[k]_{p} = J[k-1]_p$ for all $p\leq k-1$,
\item\label{Prop2} $I[k-1] \subseteq I[k]$,
\item\label{Prop3} $J[M] = I[M]$. 
\end{enumerate}
We will prove these properties respectively in Corollary \ref{CorSeqSec}, Lemma \ref{LemArgRow} and Lemma \ref{LemPrinc}.
\end{proof}

Before proceeding to the proof, we want to provide an example of possible shapes and the corresponding leading minors for the case of $N=3$. We can organize the possible shapes in terms of rank and involution $\tau$, leaving only the choice of non-zero $u$:
\begin{Exa}
Let $s_{i}\in\{1,-1\}$ and $y=e^{i\pi\theta}$ for $\theta\in[0,2\pi)$. Then up to a choice of $s_i$ and $\theta$, the possible shapes for $\mcO_q(H(3))$ are as follows:
\begin{enumerate}
	\item The rank $3$ shapes are:
	\begin{enumerate}
\item $\tau=\id$, $u([3])=\{s_1,s_2,s_3\}$, $Z_{S,1}=Z_{[1]}=Z_{11}$, $Z_{S,2}=Z_{[2]}=Z_{12,12}$, $Z_{S,3}=Z_{[3]}=Z_{123,123}$.
\item $\tau=(1,2)$, $u([3])=\{y,\bar{y},s_1\}$, $Z_{S,1}=Z_{21}$, $Z_{S,2}=Z_{[2]}$, $Z_{S,3}=Z_{[3]}$.
\item $\tau=(1,3)$, $u([3])=\{y,s_1,\bar{y}\}$, $Z_{S,1}=Z_{31}$, $Z_{S,2}=Z_{13,13}$, $Z_{S,3}=Z_{[3]}$.
\item $\tau=(2,3)$, $u([3])=\{s_1,y,\bar{y}\}$, $Z_{S,1}=Z_{32}$, $Z_{S,2}=Z_{23,23}$, $Z_{S,3}=Z_{[3]}$.
	\end{enumerate}
\item The rank $2$ shapes are:
\begin{enumerate}
	\item $\tau=\id$, $u([3])=\{s_1,s_2,0\}$, $Z_{S,1}=Z_{[1]}$, $Z_{S,2}=Z_{[2]}$.
	\item $\tau=\id$, $u([3])=\{s_1,0,s_2\}$, $Z_{S,1}=Z_{[1]}$, $Z_{S,2}=Z_{13,13}$.
	\item $\tau=\id$, $u([3])=\{0,s_1,s_2\}$, $Z_{S,1}=Z_{22}$, $Z_{S,2}=Z_{23,23}$.
	\item $\tau=(1,2)$, $u([3])=\{y,\bar{y},0\}$, $Z_{S,1}=Z_{21}$, $Z_{S,2}=Z_{[2]}$.
	\item $\tau=(1,3)$, $u([3])=\{y,0,\bar{y}\}$, $Z_{S,1}=Z_{31}$, $Z_{S,2}=Z_{13,13}$.
	\item $\tau=(2,3)$, $u([3])=\{0,y,\bar{y}\}$, $Z_{S,1}=Z_{32}$, $Z_{S,2}=Z_{23,23}$.
\end{enumerate}	
\item The rank $1$ shapes are:
\begin{enumerate}
	\item $\tau=\id$, $u([3])=\{s_1,0,0\}$, $Z_{S,1}=Z_{11}$.
	\item $\tau=\id$, $u([3])=\{0,s_1,0\}$, $Z_{S,1}=Z_{22}$.
	\item $\tau=\id$, $u([3])=\{0,0,s_1\}$, $Z_{S,1}=Z_{33}$.
\end{enumerate}
\end{enumerate}
\end{Exa} 

\begin{Rem}
In the above example, the first shape for each rank can be seen to correspond to the \emph{big cell representations} defined in \cite{DCMo24a}. From the viewpoint of shapes, a big cell representation is a direct sum of representations of shapes of the form $S = (\id,u)$ with $u = (\eta_1,\ldots,\eta_M,0,\ldots,0)$ for $\eta_i \in \{\pm 1\}$.

In the quasi-classical picture, we see that the symplectic leaf $C_u$ containing $u= \diag(\eta_1,\ldots,\eta_M,0,\ldots,0)$ will have a $T(N)$-orbit which is dense in $H^{\leq M}(N)$, the set of all hermitian matrices with rank $\leq M$. In other words, we get an injective map 
\[
\mcO(H^{\leq M}(N)) \rightarrow \mcO(T(N)),\qquad Z \mapsto T^*uT.
\]
If we want to quantize this map, the quantization of $\mcO(T(N))$ will itself depend on the precise value of $(\eta_1,\ldots,\eta_M,0,\ldots,0)$. Hence the construction of $\mcO_q^{\epsilon}(T(N))$ and the weak containment of general representations within big cell representations shown in \cite{DCMo24a} is naturally motivated by the quasi-classical picture. 
\end{Rem}

\subsection{Proof of Theorem \ref{TheoShapeIrr}}

We first prove property \eqref{Prop1} in Theorem \ref{TheoShapeIrr}.

\begin{Def}
For $M <N$ and $J \in \binom{\brN}{M}$, we let $\mathfrak{L}(J)$ be the left ideal in $\mcO_q(H(N))$ generated by 
\begin{itemize}
\item $\mathfrak{L}(J)_1 = \{Z_{I',J'}\mid I',J' \in \binom{\brN}{M} \textrm{ and }J'<J\}$, and
\item $\mathfrak{L}(J)_2 = \{Z_{I',J'}\mid I',J'\in\binom{\brN}{M+1}\textrm{ and }J \subseteq J'\}$. 
\end{itemize}
\end{Def}

Recall the map $\Phi$ in \eqref{EqDefPhi}, and write 
\[
\widetilde{\mathfrak{L}}(J) = \Phi(\mathfrak{L}(J)).
\]

\begin{Lem}
The space $\widetilde{\mathfrak{L}}(J)$ is a two-sided ideal in $\mcO_q(M_N(\C))$ with its ordinary product. 
\end{Lem}
\begin{proof}
Clearly $\Delta(\widetilde{\mathfrak{L}}(J)) \subseteq \mcO_q(M_N(\C)) \otimes \widetilde{\mathfrak{L}}(J)$, so by 
\eqref{EqRevBraid} we have that $\widetilde{\mathfrak{L}}(J)$ is a left ideal in $\mcO_q(M_N(\C))$. From \eqref{EqBraidProd} it then follows that $\widetilde{\mathfrak{L}}(J)$ is the left ideal in $\mcO_q(M_N(\C))$ generated by $\widetilde{\mathfrak{L}}(J_1) $ and $\widetilde{\mathfrak{L}}(J_2)$.

From  \eqref{EqCondR} and \eqref{EqBraidComm2} we see that $X_{I',J'}\mcO_q(M_N(\C))$ lies in the left ideal generated by all $X_{B,D}$ with $D \preceq J'$. It follows then immediately that $\widetilde{\mathfrak{L}}(J)_1 \mcO_q(M_N(\C)) \subseteq \widetilde{\mathfrak{L}}(J)$. To see that also $\widetilde{\mathfrak{L}}(J)_2 \mcO_q(M_N(\C)) \subseteq \widetilde{\mathfrak{L}}(J)$ follows from this observation, it suffices to show that $X_{B,D} \in \widetilde{\mathfrak{L}}(J)$ as soon as $D = \{d_1<\ldots<d_M<d_{M+1}\}$ with $D' = \{d_1,\ldots,d_M\} < J$. This follows from the following particular case of \eqref{EqLaplaceColumn}, where we write $B' = B \setminus \{b_{M+1}\}$,
\[
X_{B,D} = \sum_{p=1}^{M+1} (-q)^{p-M-1} X_{b_p,d_{M+1}} X_{B\setminus \{b_p\},D'}.   
\]
\end{proof}

Let us write $\mcO_q^J(M_N(\C)) = \mcO_q(M_N(\C))/\widetilde{\mathfrak{L}}(J)$ for the quotient algebra, and write the quotient map as $x \mapsto \overline{x}$. For $L \in \binom{\brN}{M}$, we write $V_L$ for the $\C$-linear space spanned by all $\overline{X}_{I,L}$ with $I \in \binom{\brN}{M}$. 

\begin{Lem}\label{LemCommVV}
We have $V_LV_J = V_JV_L$ for all  $L \in \binom{\brN}{M}$.
\end{Lem} 
\begin{proof}
This follows immediately from \eqref{EqBraidComm1} and \eqref{EqCondR} as we obtain 
\[
\overline{X}_{K,L} \overline{X}_{I,J} = \sum_{A,B,C,D} \widehat{R}^{A,I}_{K,B} (\widehat{R}^{-1})^{J,C}_{D,L} \overline{X}_{A,C} \overline{X}_{B,D} =  \sum_{A,B} \widehat{R}^{A,I}_{K,B} (\widehat{R}^{-1})^{J,J}_{L,L} \overline{X}_{A,J} \overline{X}_{B,L}.
\]
\end{proof}

We will now show that $\widetilde{\mathfrak{L}}(J)$ almost contains all quantum minors of rank $M+1$. For this, we consider again the algebra $\wedge_q(\C^N)$ with its coaction $\alpha = \wedge \delta$ defined in \eqref{EqComodAlgMapqExt}. Let us write the quotient map
\[
\overline{\alpha}: \wedge_q(\C^N) \rightarrow \wedge_q(\C^N) \otimes \mcO_q^J(M_N(\C)).
\]

\begin{Lem}\label{LemMultZero}
For any $K',L'\in \binom{\brN}{M+1}$ we have 
\[
\overline{X}_{K',L'} V_J^{M-|J \cap L'|} = 0. 
\]
\end{Lem} 
\begin{proof}
We will prove this by induction on $f = M -|J\cap L'|$. 

For $f = 0$ there is nothing to prove. Assume now that the result holds for all $f' <f$. Pick $m \in L' \setminus (J\cap L')$, and write $J' = J\cup\{m\}$. From \eqref{EqLaplaceRow} it follows that for all $1\leq s,t\leq M+1$ and $I' \in \binom{\brN}{M+1}$ we have

\begin{equation}\label{EqLinCombZero}
0 = \delta_{s,t} \overline{X}_{I',J'} = \sum_{p=1}^{M+1}(-q)^{p} \overline{X}_{i_s',j_p'} \overline{X}_{I'\setminus \{i_t'\},J'\setminus \{j_p'\}},
\end{equation}
and hence 
\begin{equation}
\sum_{p=1}^{M+1} (-q)^{p} \overline{\alpha}(e_{j_p'})_{12} \overline{\alpha}(e_{J'\setminus \{j_p'\}})_{02} = 0.
\end{equation}
In particular, 
\begin{equation}
\sum_{p=1}^{M+1} (-q)^{p} \overline{\alpha}(e_{L'\setminus \{m\}} \wedge_q e_{j_p'})_{12} \overline{\alpha}(e_{J'\setminus \{j_p'\}})_{02}(1\otimes 1 \otimes V_J^{f-1}) = 0.
\end{equation}
By induction and Lemma \ref{LemCommVV}, we know that all terms on the left hand side vanish when $j_p' \neq m$, so we deduce that 
\begin{equation}
\overline{\alpha}(e_{L'})_{12} \overline{\alpha}(e_{J})_{02}(1\otimes 1 \otimes V_J^{f-1})=0,
\end{equation}
so in particular $\overline{X}_{K',L'} V_J^f = 0$ for all $K' \in \binom{\brN}{M+1}$. 
\end{proof}

\begin{Lem}\label{LemRank} 
Let $M< N$, and let $\pi$ be an irreducible representation of $\mcO_q(H(N))$ of rank $\geq M$. Let $(J,I) \in \binom{\brN}{M}\times \binom{\brN}{M}$ be the first element with respect to $\leq$ such that $\pi(Z_{I,J}) \neq 0$. Assume moreover that $\pi(Z_{I',J'}) =0$ for all $I',J'\in   \binom{\brN}{M+1}$ with $J \subseteq J'$. Then the rank of $\pi$ is $M$. 
\end{Lem} 
\begin{proof}
Clearly $\mathfrak{L}(J) \subseteq \Ker(\pi)$, hence by Lemma \ref{LemMultZero} and \eqref{EqBraidProd} we have $Z_{K',L'} (Z_{I,J})^N \in \Ker(\pi)$ for $K',L' \in \binom{\brN}{M+1}$. As $\pi(Z_{I,J})$ is a normal operator with zero kernel and purely discrete spectrum, we deduce that $\pi(Z_{K',L'}) = 0$ for $K',L' \in \binom{\brN}{M+1}$,  hence $\pi$ has rank $M$. 
\end{proof}

The following corollary now establishes property \eqref{Prop1} in Theorem \ref{TheoShapeIrr}.

\begin{Cor}\label{CorSeqSec}
Let $\pi$ be an irreducible representation of $\mcO_q(H(N))$. Assume that the rank of $\pi$ is $\geq M+1$, and let $J' \in \binom{\brN}{M+1}$ be the $\leq$-first element with $\pi(Z_{I',J'}) \neq 0$ for some $I'\in \binom{\brN}{M+1}$. Let $J = J'\setminus \{j_{M+1}'\}$. Then $J$ is the  $\leq$-first element in $\binom{\brN}{M}$ with $\pi(Z_{I,J}) \neq 0$ for some $I\in \binom{\brN}{M}$.
\end{Cor} 
\begin{proof}
Let $\widetilde{J}$ be the  $\leq$-first element  in $\binom{\brN}{M}$ with $\pi(Z_{I,\widetilde{J}}) \neq 0$ for some $I\in \binom{\brN}{M}$. From the Laplace expansion \eqref{EqLaplExp2} for $Z_{I',J'}$ along $K = \{M+1\} \in \binom{[M+1}{1}$ we see that $\widetilde{J} \leq J$. If this equality is however strict, we have that $Z_{I',\widetilde{J}\cup \{p\}} = 0$ for all $p \in \brN \setminus \widetilde{J}$ and all $I' \in \binom{\brN}{M+1}$, so that $\pi$ has rank $M$ by Lemma \ref{LemRank}, a contradiction. Hence $\widetilde{J} = J$. 
\end{proof}

We now prove property \eqref{Prop3} in Theorem \ref{TheoShapeIrr}.

\begin{Lem}\label{LemPrinc}
Let $\pi$ be an irreducible representation of $\mcO_q(H(N))$ of rank $M$. Let $(J,I) \in \binom{\brN}{M}\times \binom{\brN}{M}$ be $\leq$-first element with $\pi(Z_{I,J}) \neq 0$. Then $I = J$.
\end{Lem} 
\begin{proof} 
Assume that $I \neq J$, so $k = |I\cup J|>M$. Denote by $T$ the symmetric difference $T = I\Delta J$ with $|T|= t$, and write $S =I \cap J$ with then $|S| = s = k-t$. Take $F \in \binom{\brk}{s}$ such that $(I\cup J)_F = I\cap J$. Finally, take $K \in \binom{\brt}{M-s}$ such that $S \cup T_K = I$, and hence also $S\cup T^K = J$. Then as $\pi$ is of rank $M$, we obtain from $L_l(I\cup J,I\cup J;F,G;K,K)$ in  \eqref{EqMuirBr} that 
\[
0 = \sum_{A,B,C,D,P} (-q)^{\wt(P) - \wt(K)} (\widehat{R}^{-1})^{A,I}_{J,B} \widehat{R}^{S\cup T_P,C}_{B,D} \pi(Z_{A,C})\pi(Z_{D,S\cup T^P}).
\] 
Recalling that $\pi(Z_{A,B})^* = \pi(Z_{B,A})$ for all $A,B$, we deduce from the assumptions that $\pi(Z_{A,B}) = \pi(Z_{B,A}) =0$ for $B\prec J$, and hence we obtain from \eqref{EqCondR} that 
\[
0  =  \sum_{P} (-q)^{\wt(P) - \wt(K)} (\widehat{R}^{-1})^{I,I}_{J,J} \widehat{R}^{S\cup T_P,S\cup T_P}_{J,J} \pi(Z_{I,S\cup T_P})\pi(Z_{J,S\cup T^P}).
\]
Now the sum need only be taken over those $P \in  \binom{\brt}{M-s}$ with $S \cup T_P \geq J$ and $S\cup T^P \geq I$. By the next lemma, this means the only non-trivial term arises when $S\cup T_P = J$ and hence $S\cup T^P = I$. We arrive at 
\[
0  =  \pi(Z_{I,J})\pi(Z_{J,I}) = \pi(Z_{I,J})\pi(Z_{I,J})^*.
\]
Hence $\pi(Z_{I,J})=0$, a contradiction.
\end{proof}

\begin{Lem}\label{LemCombPart}
Assume $I \in \binom{\brN}{M}$ and $J \in \binom{\brN}{M'}$. Let $S = I\cap J$ and $T = I \Delta J$ with $|S|=s$ and $|T|=t$. Assume $P \in \binom{\brt}{M'-s}$ with $S\cup T_P \geq J$ and $S\cup T^P \geq I$. Then both inequalities must be equalities. 
\end{Lem}
\begin{proof}
We have 
\[
S \cup T_P \geq J= S\cup (J\setminus S),
\] 
hence $T_P \geq J\setminus S$. Similarly $T^P \geq I\setminus S$, so it is enough to prove the lemma assuming $S = \emptyset$. Then  $P \in \binom{[M+M']}{M'}$ with $T_P \geq J$ and $T^P \geq I$. We write in the following $P' = T \setminus P$, so $T^P = T_{P'}$. 

We prove this by induction, taking as induction basis the trivial case $I= \emptyset$ or $J = \emptyset$. Assuming $I,J$ both non-trivial, we may by symmetry suppose $j_1< i_1$, so that necessarily $t_{p_1} = j_1$.  

We claim that also $t_{p_1'} = i_1$. Indeed, suppose first that $t_{p_1'} = j_p$ for some $p>1$. As $I \leq T^P$ we then have $i_1 < j_p$, so necessarily $i_1 \in T_P$. Since there are then also only $p$ elements $\leq i_1$ in $I \cup J$, we must have $t_{p_a} = i_1$ with $a \leq p$. If $a<p$ we have 
\[
T_P = \{j_1<j_2<\ldots < j_{a-1} < i_1 < \ldots < j_a < \ldots\},
\]
which by $J \leq T_P$ leads to the contradiction $j_a \leq i_1 < j_a$. If $a=p$ we have 
\[
T_P = \{j_1<j_2<\ldots < j_{p-1} < i_1 < \ldots\},
\]   
which again by $J \leq T_P$ leads to the contradiction $j_p \leq i_1 <j_p$. 

In conclusion, we find that we must have $t_{p_1'} = i_p$ for some $p$. If $p \neq 1$, then again $i_1 \in T_P$. We must then have 
\[
T_P = \{j_1<\ldots < j_a < i_1<\ldots < j_{a+1} < \ldots\},
\]
where the first string must contain all $j_b$ with $1\leq b \leq a$ since otherwise we get the contradiction $i_p < j_b <i_1 \leq i_p$. But from $J\leq T_P$ we see then again that $j_{a+1} \leq i_1 < j_{a+1}$, a contradiction. Hence $p=1$. 

We have obtained now $T_P = \{j_1< \ldots\}$ and $T^P = \{i_1<\ldots\}$. Hence also $T_P \setminus \{j_1\} \geq J \setminus \{j_1\}$ and $T^P \setminus \{i_1\} \geq I \setminus \{i_1\}$. By induction, we are done. 
\end{proof}

Finally, we prove property \eqref{Prop2} in Theorem \ref{TheoShapeIrr}.

\begin{Lem}\label{LemArgRow}
Let $\pi$ be an irreducible representation of $\mcO_q(H(N))$ of rank $>M$, and let $J' \in \binom{\brN}{M+1}$ be the $\leq$-first element with $\pi(Z_{I',J'}) \neq 0$ for some $I'$. Let $J = J' \setminus \{j_{M+1}\}$. Let $I$ be the $\leq$-first element with $\pi(Z_{I,J}) \neq 0$, and let $I'$ be the $\leq$-first element with $\pi(Z_{I',J'})\neq 0$. Then $I \subset I'$.  
\end{Lem}
\begin{proof}
We prove this by contradiction. Assume $I \not\subset I'$, and put $k = |I\cup I'|$ and $r = |I\Delta I'|$. Put $0 < l = |I \setminus (I\cap I')|$. By cardinality considerations we may then pick an arbitrary $R \subseteq \brN$ disjoint from $J'$ such that $\widetilde{J} := J' \cup R$ has $|\widetilde{J}| = k$, as well as $F,G\in \binom{\brk}{k-r}$ and $K,K'\in \binom{r}{l}$ such that 
\begin{itemize}
\item $(I \cup I')_F = I \cap I'$, 
\item $\widetilde{J}_G \cup(\widetilde{J}^{G})_K =J$, and
\item $\widetilde{J}_G \cup (\widetilde{J}^G)^{K'} = J'$.
\end{itemize}
Then using \eqref{EqMuirBr2} in the form $L_r(I\cup I',\widetilde{J};F,G;K,K')$ we find
\[
0 = \sum_{A,B,C,D} \sum_{P \in \binom{\brr}{l}} (-q)^{\wt(P) -\wt(K)} (\widehat{R}^{-1})^{A,(I\cap I')\cup (I\Delta I')_P}_{(I\cap I')\cup (I \Delta I')^{P},B} \widehat{R}^{J,C}_{B,D} \pi(Z_{A,C}) \pi(Z_{D,J'})
\]
which by the assumptions and the properties of $\widehat{R}$ simplifies to 
\begin{equation}\label{EqZeroExp}
0 = \sum_{P \in \binom{\brr}{l}}\sum_{A,D}  (-q)^{\wt(P) -\wt(K)} (\widehat{R}^{-1})^{A,(I\cap I')\cup (I\Delta I')_P}_{(I\cap I')\cup (I \Delta I')^{P},D} \widehat{R}^{J,J}_{D,D} \pi(Z_{A,J}) \pi(Z_{D,J'}).
\end{equation}
Let us now consider a term for fixed $P$, and write for simplicity (recycling notation) $F = I \cap I'$, $K = (I\Delta I')_P$ and $L = (I \Delta I')^P$. Then by the properties of the $R$-matrix, we only need to sum over those $A,D$ with 
\[
(F\cup K) \setminus A = D \setminus (F\cup L), \qquad A \setminus (F\cup K) = (F\cup L) \setminus D.\]
This clearly implies that $A \cap D = I\cap I'$ and $A\cup D = I \cup I'$, so we can write $A = (I\cap I')\cup (I\Delta I')_{P'}$ and $D= (I\cap I')\cup (I\Delta I')^{P'}$ for some $P'$. But we may in the summation also restrict to terms $A,D$ with
\[
I \leq A,\qquad I' \leq D,
\]
in which case $A=I$ and $D=I'$ by Lemma \ref{LemCombPart}. 

We can now simplify \eqref{EqZeroExp} to 
\[
0 = \left(\sum_{P \in \binom{\brr}{l}} (-q)^{\wt(P)} (\widehat{R}^{-1})^{I,(I\cap I')\cup (I\Delta I')_P}_{(I\cap I')\cup (I \Delta I')^{P},I'}\right) \widehat{R}^{J,J}_{I',I'} \pi(Z_{I,J}) \pi(Z_{I',J'}).
\]
As $\pi(Z_{I,J})$ and $\pi(Z_{I',J'})$ are non-zero $q$-commuting normal elements with zero kernel and purely discrete spectrum, we will obtain a contradiction once we know that 
\[
\sum_{P \in \binom{\brr}{l}} (-q)^{\wt(P)} (\widehat{R}^{-1})^{I,(I\cap I')\cup (I\Delta I')_P}_{(I\cap I')\cup (I \Delta I')^{P},I'}\neq 0.
\]
But 
\begin{multline}\label{EqRNonZero}
\sum_{P \in \binom{\brr}{l}} (-q)^{\wt(P)} (\widehat{R}^{-1})^{I,(I\cap I')\cup (I\Delta I')_P}_{(I\cap I')\cup (I \Delta I')^{P},I'} 
\\= \sum_{P \in \binom{\brr}{l}} (-q)^{\wt(P)} \langle e_{(I\cap I')\cup (I\Delta I')_P} \otimes e_{(I\cap I')\cup (I\Delta I')^P},\widehat{\msR}^{-1}(e_{I'} \otimes e_{I})\rangle. 
\end{multline}
which is non-zero by self-adjointness of $\widehat{\msR}$ and the following Lemma \ref{LemCompRmatrix} beneath.
\end{proof}

Here we view $\wedge_q^k(\C^N)=\{e_I\}$ as an irreducible $U_q(\mfgl(N,\C))$ module via the dual pairing on $\mcO_q(U(N))$ \cite{DCMo24a}. From this viewpoint, the skew bicharacter $\mbr$ becomes a universal $R$-matrix in the completion of $U_q(\mfgl(N,\C))\otimes U_q(\mfgl(N,\C))$, which we denote by $\widehat{\msR}$.

\begin{Lem}\label{LemCompRmatrix}
Let $I,I' \subseteq [N]$ with $|I| = M$, $|I'|= M'$, and further write 
\[
l = |I \setminus I'|,\qquad l' = |I'\setminus I|,\qquad r = |I \Delta I'| = l+l',
\] 
where $\Delta$ denotes the symmetric difference. Put 
\[
\xi = \sum_{P \in \binom{\brr}{l}} (-q)^{\wt(P)} e_{(I\cap I')\cup (I\Delta I')_P} \otimes e_{(I\cap I')\cup (I\Delta I')^P} \in \wedge_q^M(\C^N) \otimes \wedge_q^{M'}(\C^N),
\]
\[
\xi'= \sum_{P' \in \binom{\brr}{l'}} (-q)^{\wt(P')} e_{(I\cap I')\cup (I\Delta I')_{P'}} \otimes e_{(I\cap I')\cup (I\Delta I')^{P'}} \in \wedge_q^{M'}(\C^N) \otimes \wedge_q^{M}(\C^N).
\]
Then
\[
\widehat{\msR}^{-1}\xi =   q^{|I\cap I'|}(-q)^{\frac{l(l+1)}{2} - \frac{l'(l'+1)}{2}- ll'}\xi'.
\]
\end{Lem}
\begin{proof}
As the category of integral $U_q(\mfgl(N,\C))$-modules is braided by the $\widehat{\msR}_{V,W}$, we have the following $U_q(\mfgl(N,\C))$-module algebra structures on the algebraic tensor product $\wedge_q(\C^N) \otimes \wedge_q(\C^N)$: 
\[
(a\otimes c) \cdot (b\otimes d) = (a\wedge_q(\msR_{2}b)) \otimes ((\msR_1 c)\wedge_q d),
\]
\[
(a\otimes c) * (b\otimes d) = (a\wedge_q(\msR_{1}^{-1} b)) \otimes ((\msR_2^{-1} c)\wedge_q d),
\]
where we use the shorthand `$\msR = \msR_1\otimes \msR_2$'. Moreover, an easy diagrammatic computation, using the naturality of $\widehat{\msR}$, gives that 
\begin{equation}\label{EqBraidProdSwitch}
\widehat{\msR}^{-1}((a\otimes b) \cdot (c\otimes d)) = (\widehat{\msR}^{-1}(a\otimes b))* (\widehat{\msR}^{-1}(c\otimes d)). 
\end{equation}

Now for $A \cap B = \emptyset$ we have that $\msR (e_A \otimes e_B) - e_A \otimes e_B$ is a linear combination of $e_{A'}\otimes e_{B'}$ with $A' \sqcup B' = A \sqcup B$ and $A' \neq A$ and $B' \neq B$. It follows that
\begin{eqnarray} 
&& \hspace{-2cm} (e_{I\cap I'} \otimes e_{I\cap I'}) \cdot (\sum_{P \in \binom{\brr}{l}} (-q)^{\wt(P)} e_{(I\Delta I')_P} \otimes e_{(I\Delta I')^P})  \nonumber\\
&=&  \sum_{P \in \binom{\brr}{l}} (-q)^{\wt(P)}(e_{I\cap I'} \wedge_q (\msR_2 e_{(I\Delta I')_P})) \otimes ((\msR_1e_{I\cap I'})\wedge_q e_{(I\Delta I')^P}) \nonumber\\ &= &   \sum_{P \in \binom{\brr}{l}} (-q)^{\wt(P)}(e_{I\cap I'} \wedge_q e_{(I\Delta I')_P}) \otimes (e_{I\cap I'}\wedge_q e_{(I\Delta I')^P}) \nonumber\\ 
 &= & (-q)^{\wt(K) - \frac{s(s+1)}{2}}  \xi \label{EqIdProdDot},
\end{eqnarray} 
where $s =  |I \cap I'|$ and $K \in \binom{[r+s]}{s}$ is such that $(I\cup I')_K = I\cap I'$.

By \eqref{EqBraidProdSwitch}, it follows that 
\begin{equation}\label{EqCompRr}
\widehat{\msR}^{-1}\xi = (-q)^{-\wt(K) + \frac{s(s+1)}{2}} (\widehat{\msR}^{-1}(e_{I\cap I'} \otimes e_{I\cap I'})) *( \widehat{\msR}^{-1} \sum_{P \in \binom{\brr}{l}} (-q)^{\wt(P)} e_{(I\Delta I')_P} \otimes e_{(I\Delta I')^P}).
\end{equation}
Now it is easily seen that $\widehat{\msR}^{-1}(e_{I\cap I'} \otimes e_{I\cap I'}) = q^{|I\cap I'|}e_{I\cap I'} \otimes e_{I\cap I'}$. On the other hand, we claim that 

\begin{equation}\label{EqRhatActs}
 \widehat{\msR}^{-1} \sum_{P \in \binom{\brr}{l}} (-q)^{\wt(P)} e_{(I\Delta I')_P} \otimes e_{(I\Delta I')^P} = (-q)^{\frac{l(l+1)}{2} - \frac{l'(l'+1)}{2}-ll'}  \sum_{P' \in \binom{\brr}{l'}} (-q)^{\wt(P')} e_{(I\Delta I')_{P'}} \otimes e_{(I\Delta I')^{P'}}.
\end{equation}
Indeed, the projection $\rho_l:  (\C^r)^{\otimes l} \rightarrow \wedge_q^l(\C^r)$ has a splitting as an $U_q(\mfgl(r,\C))$-module via
\[
\iota_l: e_{P} \mapsto \frac{1}{[l]_{q^2}!} \sum_{\sigma \in \Sym(P)} (-q)^{l(\sigma)}e_{\sigma(p_1)}\otimes \ldots \otimes e_{\sigma(p_l)},\qquad P = \{p_1<\ldots<p_l\}
\]
where $[l]_{q^2}! = \prod_{k=1}^l \frac{(1-q^{2k})}{(1-q^2)}$. Hence, with $P^c = \{p_1'<\ldots<p_{l'}'\}$ denoting the complement of $P$ in $[r]$ and with $\sigma_P$ the permutation
\[
\sigma_P: (1,2,\ldots,r) \mapsto (p_1,\ldots,p_l,p_1',\ldots,p_{l'}')
\] 
of length $\wt(P) -l(l+1)/2$, we see that
\begin{eqnarray*}
&& \hspace{-1cm} [l]_{q^2}![l']_{q^2}!\sum_{P \in \binom{\brr}{l}} (-q)^{\wt(P)} \iota_l(e_{P}) \otimes \iota_{l'}(e_{P^c}) \\ &= &  \sum_{P \in \binom{\brr}{l}}  \underset{\sigma'\in \Sym(P^c)}{\sum_{\sigma \in \Sym(P)}} (-q)^{\wt(P) + l(\sigma)+l(\sigma')} (e_{\sigma(p_1)}\otimes \ldots\otimes e_{\sigma(p_l)})\otimes (e_{\sigma'(p_1')}\otimes \ldots \otimes e_{\sigma'(p_{l'}')})\\ &=&\sum_{P \in \binom{\brr}{l}} \underset{\sigma'\in \Sym(P^c)}{\sum_{\sigma \in \Sym(P)}} (-q)^{\wt(P)+l(\sigma)+l(\sigma')} (e_{\sigma\sigma_P(1)}\otimes \ldots \otimes e_{\sigma\sigma_P(l)})\otimes (e_{\sigma'\sigma_P(l+1)}\otimes \ldots \otimes e_{\sigma'\sigma_P(r)}) \\
&=& (-q)^{\frac{l(l+1)}{2}}\sum_{\sigma \in \Sym([r])} (-q)^{l(\sigma)} (e_{\sigma(1)}\otimes \ldots\otimes e_{\sigma(l)})\otimes (e_{\sigma(l+1)}\otimes \ldots \otimes e_{\sigma(r)}),
\end{eqnarray*}
which is a $-q$-eigenvector under all $\hat{R}_{ij}$. So
\begin{eqnarray*}
&& \hspace{-1cm}\widehat{\msR}^{-1} \sum_{P \in \binom{\brr}{l}} (-q)^{\wt(P)} e_{P} \otimes e_{P^c}\\ &=& (\rho_{l'}\otimes \rho_l)\left(\widehat{\msR}^{-1}\sum_{P \in \binom{\brr}{l}} (-q)^{\wt(P)} \iota_l(e_{P}) \otimes \iota_{l'}(e_{P^c})\right)\\ &=& (-q)^{-ll'} (\rho_{l'}\otimes \rho_l)\left(\frac{(-q)^{l(l+1)/2}}{[l]_{q^2}![l']_{q^2}!}\sum_{\sigma \in \Sym([r])} (-q)^{l(\sigma)} (e_{\sigma(1)}\otimes \ldots\otimes e_{\sigma(l')})\otimes (e_{\sigma(l'+1)}\otimes \ldots \otimes e_{\sigma(r)})\right)\\ &=& (-q)^{\frac{l(l+1)}{2} - \frac{l'(l'+1)}{2}-ll'}  \sum_{P' \in \binom{\brr}{l'}} (-q)^{\wt(P')} e_{P'} \otimes e_{(P')^c}.
\end{eqnarray*}
This proves \eqref{EqRhatActs} (using that maps of the form $e_{i_1}\otimes \cdots e_{i_l} \mapsto e_{\kappa(i_1)} \otimes \cdots \otimes e_{\kappa(i_l)}$, for $\kappa: [r] \rightarrow [N]$ an order preserving injection, intertwine the action of the braid operator for $U_q(\mfgl(r,\C))$ and the braid operator for $U_q(\mfgl(N,\C))$).

Plugging the above identities into \eqref{EqCompRr} we find, by a similar computation as for \eqref{EqIdProdDot} that 
\begin{eqnarray*}
&&\hspace{-1cm} \widehat{\msR}^{-1}\xi  \\ &=& q^{|I\cap I'|}(-q)^{-\wt(K)+\frac{s(s+1)}{2} +\frac{l(l+1)}{2} - \frac{l'(l'+1)}{2}- ll'} (e_{I\cap I'} \otimes e_{I\cap I'}) *(\sum_{P' \in \binom{\brr}{l'}} (-q)^{\wt(P')} e_{(I\Delta I')_{P'}} \otimes e_{(I\Delta I')^{P'}})\\
&=& q^{|I\cap I'|}(-q)^{\frac{l(l+1)}{2} - \frac{l'(l'+1)}{2}- ll'} \xi'.
\end{eqnarray*}
\end{proof}

\end{document}